\def\draft{n}
\documentclass{amsart}
\usepackage{fullpage,amssymb,epic,eepic,epsfig,amsmath}

\theoremstyle{plain}

\newtheorem{theorem}{Theorem}
\newtheorem{proposition}{Proposition}[section]
\newtheorem{lemma}[proposition]{Lemma}
\newtheorem{corollary}[proposition]{Corollary}

\theoremstyle{definition}
\newtheorem{definition}[proposition]{Definition}

\theoremstyle{remark}

\newtheorem{remark}[proposition]{Remark}

\def\printname#1{
	\if\draft y
		\smash{\makebox[0pt]{\hspace{-0.5in}
			\raisebox{8pt}{\tt\tiny #1}}}
	\fi
}

\newcommand{\psdraw}[2]
         {\begin{array}{c} \hspace{-1.3mm}
	\raisebox{-4pt}{\epsfig{figure=draws/#1.eps,width=#2}}
	\hspace{-1.9mm}\end{array}}

\newlength{\standardunitlength}
\catcode`\@=11
\long\def\@makecaption#1#2{%
     \vskip 10pt

\setbox\@tempboxa\hbox{
       \small\sf{\bfcaptionfont #1. }\ignorespaces #2}%
     \ifdim \wd\@tempboxa >\captionwidth {%
         \rightskip=\@captionmargin\leftskip=\@captionmargin
         \unhbox\@tempboxa\par}%
       \else
         \hbox to\hsize{\hfil\box\@tempboxa\hfil}%
     \fi}
\font\bfcaptionfont=cmssbx10 scaled \magstephalf
\newdimen\@captionmargin\@captionmargin=2\parindent
\newdimen\captionwidth\captionwidth=\hsize
\catcode`\@=12

\def\lbl#1{\label{#1}\printname{#1}}

\def\BE{\mathbb E}

\def\s{\sigma}

\def\w{\omega}

\def\d{\delta}

\def\s{\sigma}

\def\sgn{\operatorname{sgn}}

\def\Sym{\mathrm{Sym}}

\def\Per{\mathrm{per}}

\def\rot{\mathrm{rot}}

\def\sign{\mathrm{sign}}

\def\bq{\bar{q}}

\def\unknot{\mathrm{unknot}}

\begin{document}


\title[A permanent formula for the Jones polynomial]{A permanent formula for
the Jones polynomial}

\author{Martin Loebl}
\address{Dept.~of Applied Mathematics and\\
Institute of Theoretical Computer Science (ITI)\\
Charles University \\
Malostranske n. 25 \\
118 00 Praha 1 \\
Czech Republic.}
\email{loebl@kam.mff.cuni.cz}
\author{Iain Moffatt}
\address{Department of Mathematics and Statistics\\
University of South Alabama\\
Mobile\\
AL 36688 \\
USA.}
\email{imoffatt@jaguar1.usouthal.edu}

\thanks{M.L. gratefully acknowledges the support of CONICYT via Anillo en Redes ACT08.
\newline
2010 {\em Mathematics Classification.} Primary 57M15. Secondary 57M27, 05C10.
}

\date{ \today}

\begin{abstract}
The permanent of a square matrix is defined in a way similar to the 
determinant, but without using signs. The exact computation of the permanent
is hard, but there are Monte-Carlo algorithms that can estimate general
permanents. 
Given a planar diagram of a link $L$ with
$n$ crossings, we define a $7n \times 7n$ matrix whose permanent equals
to the Jones polynomial of $L$. This result accompanied with recent work
of Freedman, Kitaev, Larson and Wang  provides a Monte-Carlo
 algorithm to any decision problem belonging to the class BQP, i.e. such that
it can be computed with bounded error in polynomial time using quantum resources.
\end{abstract}

\maketitle



\section{Introduction and statement of results}
\lbl{sec.intro}

The {\em permanent} of an $n \times n$ matrix $A=(a_{ij})$ is defined 
to be
$$
\Per(A)=\sum_{\s \in \Sym_n} \prod_{i=1}^n a_{i \s(i)}
$$
where $\Sym_n$ is the permutation group on $\{1,\dots,n\}$. 
 It is well known that if $A$ is the $(V\times W)$, $0,1$
adjacency matrix of a bipartite graph $G= (V,W,E)$, then $\Per(A)$ is the number of 
perfect matchings of $G$.
The permanent of $A$
is syntactically similar to its determinant $\det(A)$, which is
a signed variation of the above sum. This mild sign variation leads
to a radical change in computability: 
computing permanents is hard (see \ref{sub.ex}), whereas determinants can be computed
in a polynomial time.

It is a seminal result of Valiant (see \cite{V3}) that many graph and knot polynomials, including the Jones polynomial, 
may be written as permanents. 
Here we are interested in the expression of  the Jones polynomial as a permanent and in the implications of this expression. In the case of the Jones polynomial, the general reduction of \cite{V3} leads to matrices of size  at least $n^2 \times n^2$, $n$ being the number of crossings of the link diagram.  The large size of these matrices severely restricts a computational applicability of this result. It is clear that in order to efficiently calculate the Jones polynomial as a permanent, a new approach is needed. In this paper we exploit a well-known combinatorial interpretation of the permanent to find an expression of the Jones polynomial as the permanent of a matrix that grows linearly in the number of crossings.

\medskip

 The {\em Jones polynomial} (\cite{J1}) is a  
 celebrated invariant of links in $S^3$.  A {\em link} is a disjoint union of embedded circles in $3$-space.
The Jones polynomial $J$ of a link can be uniquely characterized by the 
following
skein relation:
$$
q^2 J\left(\psdraw{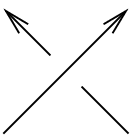}{.3in}\right)-q^{-2} J\left(\psdraw{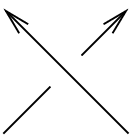}{.3in}\right)
=(q^{}-q^{-1}) J\left(\psdraw{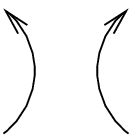}{.3in}\right)
$$
together with the initial condition $J(\unknot)=q^{}+q^{-1}$. 
From this skein relation, it follows that the Jones polynomial of a link can be computed in exponential
time (with respect to the number of crossings). Although the skein relation above provides the best known definition of the Jones polynomial, there are several other ways to construct the Jones polynomial. Below we will use a statistical mechanical construction of the Jones polynomial that is due to Turaev \cite{Tu} and Jones \cite{J2}. This state sum formulation for the Jones polynomial is described in Subsection~\ref{sub.state}.

In this paper we provide a permanent formula for the Jones polynomial and discuss applications and implications of this formula. This formulation of the Jones polynomial as a permanent  of a $7n\times 7n$ matrix (where $n$ is the number of crossings of the link diagram) is described below.

Consider a diagram $D_L$ of an oriented link $L$, that is an  oriented, $4$-valent plane graph, where each vertex has a  crossing structure of one of  two types:
\[ \text{positive: }\psdraw{L+}{0.4in}; \quad \quad \text{or negative: }\psdraw{L-}{0.4in}.\]
We form a graph $\hat{D}_L$,   from the link diagram $D_L$, by replacing a neighbourhood of each  crossing of $D_L$  with one of the  graphs shown in Figure~\ref{fig:gadget}.  
We say that $\hat{D}_L$ is a blown-up version of $D_L$. $\hat{D}_L$ is an immersed 
directed graph with $7n$ vertices, where $n$ is the number of crossings of $D_L$. 
We will refer to the the two graphs shown in Figure~\ref{fig:gadget}, minus the incoming and outgoing edges (which come from the link diagram), as {\em gadgets}. 
We emphasize the fact that the four edges in $\hat{D}_L$ that enter and exit a gadget are all parallel when they meet the gadget.

In Definition \ref{def.ww} below, we define local 
weights  on the edges of $\hat{D}_L$. Let $M_L$
denote the adjacency matrix of weights of $\hat{D}_L$. $M_L$ has size the number of
vertices of $\hat{D}_L$, and the $(i,j)$ entry of
$M_L$ is the weight of the corresponding directed edge $(ij)$ of $\hat{D}_L$.

\begin{figure}[h]
\begin{center}
\raisebox{1.5cm}{\epsfig{file=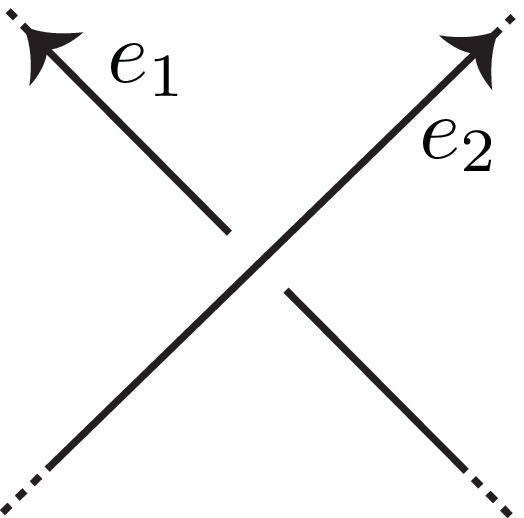, height=2cm}}\hspace{1cm} \raisebox{2.5cm}{\epsfig{file=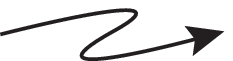, width=15mm}}  \hspace{1cm}  \epsfig{file=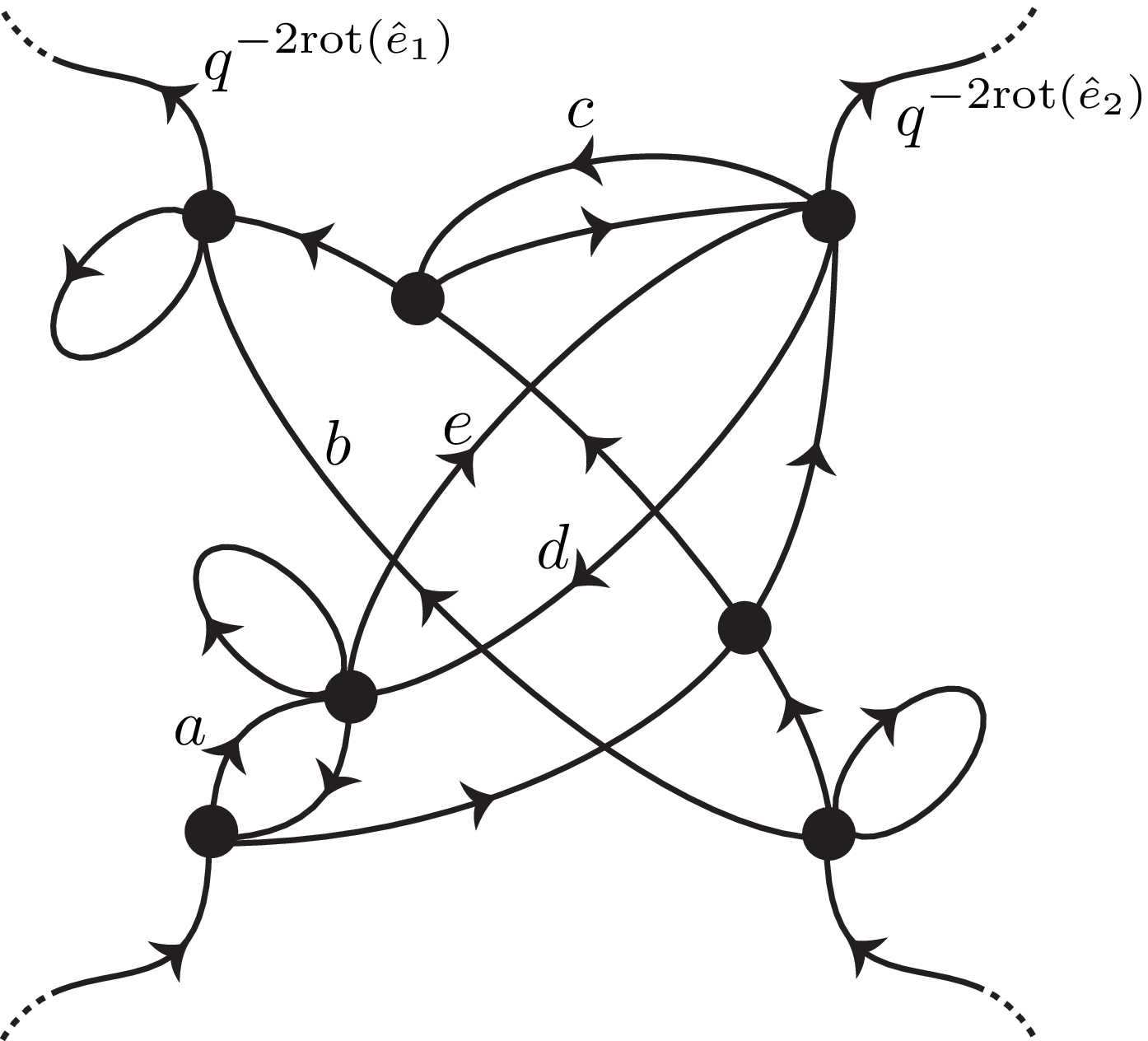, height=5cm} 
\end{center}
\vspace{1cm}
\begin{center}
\raisebox{1.5cm}{\epsfig{file=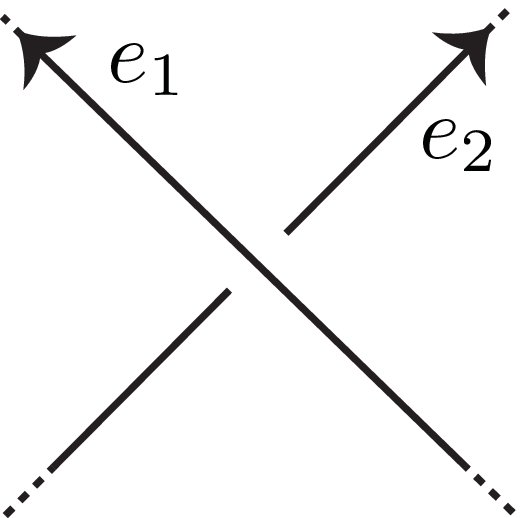, height=2cm}}\hspace{1cm} \raisebox{2.5cm}{\epsfig{file=arrow, width=15mm}}  \hspace{1cm}  \epsfig{file=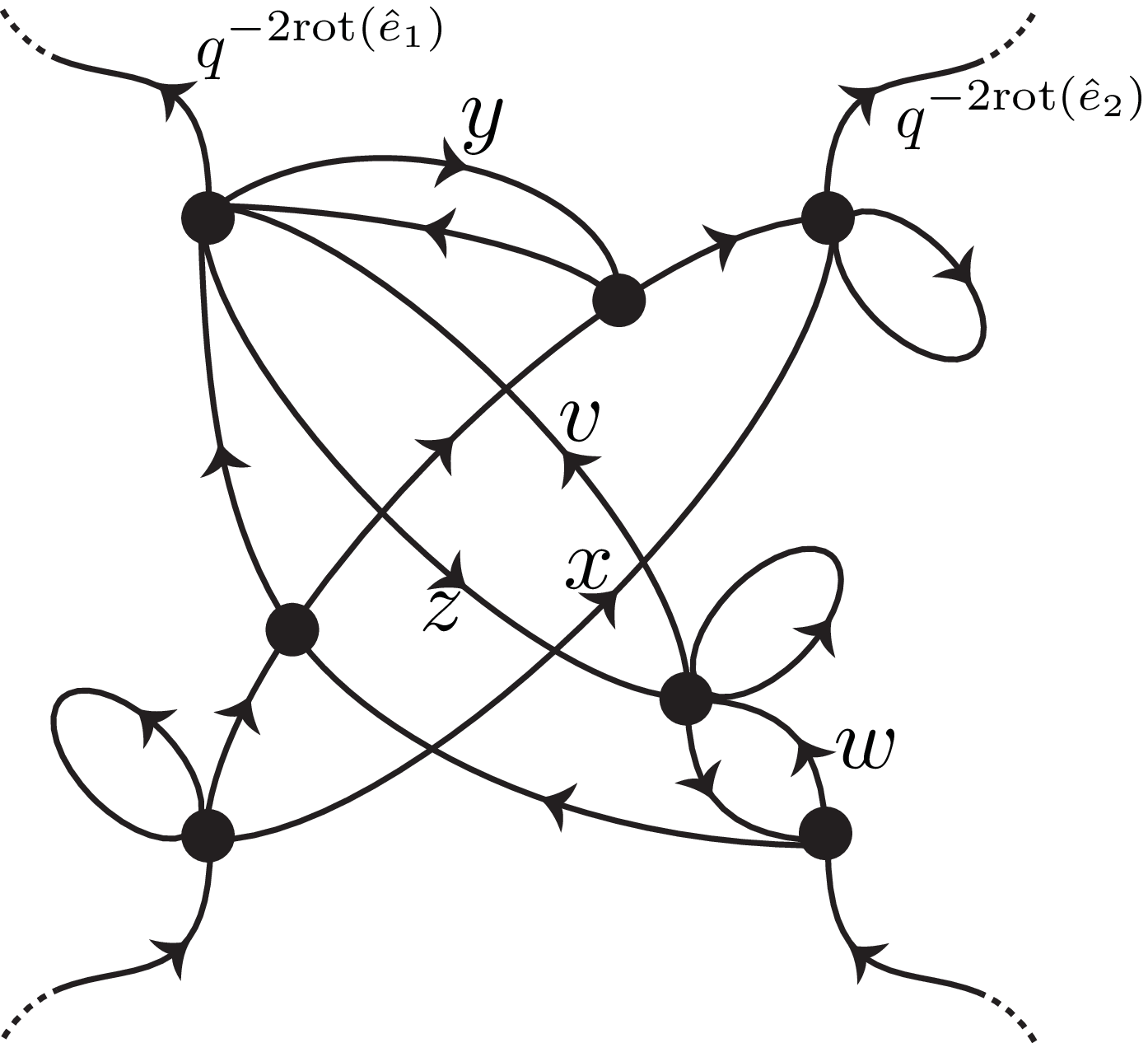, height=5cm} 
\end{center}
\caption{Gadgets at  positive and negative crossings.}
\lbl{fig:gadget}
\end{figure}

The next theorem is the main result of the paper. This result was inspired by the observation (see \cite{BG}, \cite{BL}) 
that the weight system associated with the {\em colored Jones function} is a permanent.

\begin{theorem}
\lbl{thm.per}
For every link $L$ we have
$$
J(L)= q^{-2 \w(D_L)} q^{\rot(L)}\Per(M_L).
$$
\end{theorem}
Definitions of the rotation number $\rot(L)$ and writhe $\w(D_L)$ are given in Subsection~\ref{sub.state}.

\subsection*{Acknowledgement}
The first author wants to thank D. Aharonov, A. Barvinok, S. Basu, P. Buergisser, S. Garoufalidis, M. Jerrum, M. Kiwi and P. Tetali
for enlightening discussions.

\section{Computational implications}

\subsection{Jones and quantum computing}
\lbl{sub.jqc}
We follow an exposition of the results of Freedman, Kitaev, Larsen and Wang (see \cite{FKLW}) written  by 
Bordewich, Freedman, Lovasz and Welsh (see in particular the proof of Theorem 5.1 in \cite{BFLW}). 

Suppose that we have a BQP language and an input $x$. Then we can construct a link $L$,
of size polynomial in $|x|$, such that if $x$ is in the language then
$|J(L, e^{2\pi i/5})|< [2]_5^{|x|+1}0.39$. On the other hand, if $x$ is not in the language then $|J(L, e^{2\pi i/5})|> [2]_5^{|x|+1}0.65$. We use the standard notation
$[2]_5= 2\cos (\pi/5)$.

It is also shown in \cite{FKLW} that this kind of approximation of the Jones polynomial
(called the {\em additive approximation} in \cite{BFLW}) is in BQP.
This was proven in a different way also by Aharonov, Jones and Landau \cite{AJL}. 
The BQP-hardness of approximating the Jones polynomial is also discussed by Aharonov and Arad in \cite{AA}.

\subsection{Exact computation of Jones polynomial and permanent}
\lbl{sub.ex}
The complexity class $\#P$ consists of the counting versions of the decision problems in $NP$; an example
of a problem in $\#P$ is: given a graph, how many Hamiltonian cycles does it have? It is considered
very unlikely that $\#P= P$.
Exact computation of the permanent of $0,1$ matrices is \#P-complete (see Valiant \cite{V1})
and exact computation of the Jones polynomial $J(L,t)$ is $\#P-$hard except when $t$ is a root of unity of order $r\in \{1,2,3,4,6\}$ (see Jaeger, Vertigan and Welsh \cite{JVW}). 
Kuperberg (see \cite{Kup})
showed that stronger than additive approximations of the Jones polynomial are hard. This contrasts with
the approximations of the permanents, as we explain next.

\subsection{Monte-Carlo algorithms for a permanent}
\lbl{sub.MCa}
Jerrum, Sinclair and Vigoda constructed {\em fully polynomial randomized 
approximation scheme} (FPRAS, in short) for approximating permanents of 
matrices with {\em nonnegative entries}, \cite{JSV}. For an introduction to
FPRAS, the reader may consult \cite[p. xvii]{H}, \cite{B} and \cite{JSV}.
Unfortunately, the result of
\cite{JSV} is the best possible. If one could approximate by a FPRAS
the permanent of matrices with (say) integer entries, then $P=\#P$.

There are however several approximation algorithms for general permanents, which do not 
yield a polynomial-time complexity bound. Perhaps the simplest one is based on the following
theorem, noticed by several researchers (Hammersley, Heilman, Lieb, Gutman, Godsil).

\begin{theorem}
\lbl{thm.est}
Let $A$ be a matrix and let $B$ be the random matrix obtained from $A$ by
taking the square root of minimal argument of each 
non-zero entry and then
multiplying each non-zero entry by an element of $\{1,-1\}$ chosen
independently uniformly at random. Then $\BE((\det(B))^2)= \Per(A)$.
\end{theorem}

This leads to a Monte-Carlo algorithm for estimating the permanent.
The algorithm was described and studied first by Karmarkar, Karp, Lipton, Lovasz
and Luby (see \cite{KKLLL}). 
Clearly, for general matrices the Monte-Carlo algorithm described above may
have to run an exponential time; however, with regard to the connection to the 
quantum computing, the experimental study of the algorithm applied to
the particular matrices which come from knot and link diagrams is
an attractive task and it is our work in progress jointly with Petr Plechac.

\section{Proofs}
\lbl{sec.proofs}

\subsection{State sums and gadgets}
\lbl{sub.state}

The Jones polynomial has deep connections with statistical mechanics (see \cite{J2,Wu} for  example). In particular, the Jones polynomial can be defined as  an `ice-type' or `vertex' statistical mechanical model. (The Jones polynomial can also be defined using other types of statistical mechanical models.) This means that the  Jones polynomial can be defined 
 as a state sum
which uses  a planar diagram of a link. This ice-type model for the Jones polynomial is due to Turaev \cite{Tu} and Jones \cite{J2} and is described below. We will use this model in our construction of a permanent formula of the Jones polynomial.

Choose a diagram $D_L$ of a link $L$.
Then, the Jones polynomial is given by the state sum
$$
J(L)(q)=q^{-2 \w(D_L)} \sum_{s} 
q^{\rot_0(s)-\rot_1(s)} \prod_v R^{\sgn(v)}_v(s).
$$  
where
\begin{itemize}
\item
A {\em state} $s$ is a coloring of the edges of $D_L$ by $0$ or $1$ such that
around each vertex (positive or negative) of $D_L$ the coloring looks like 
one of the following possibilities:
$$
\psdraw{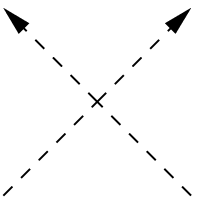}{0.4in} \hspace{0.5cm}
\psdraw{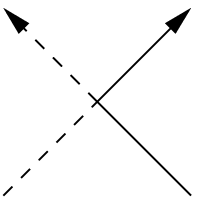}{0.4in} \hspace{0.5cm}
\psdraw{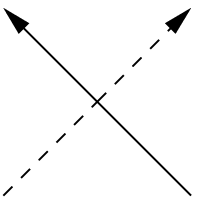}{0.4in} \hspace{0.5cm}
\psdraw{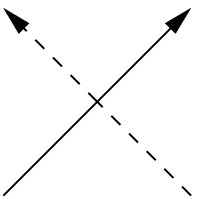}{0.4in} \hspace{0.5cm}
\psdraw{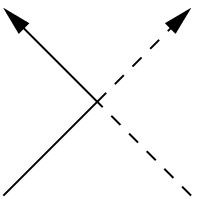}{0.4in} \hspace{0.5cm}
\psdraw{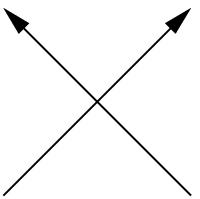}{0.4in} \hspace{0.5cm}
$$
where edges colored by $0$ or $1$ are depicted as dashed or solid 
respectively.
\item
The {\em local weight} $R^{\pm}_v(s)$ of a state $s$ at a vertex $v$
of $D_L$ (positive or negative) is given by
\begin{gather}
\text{
\begin{tabular}{|c|c|c|c|c|c|c|} \hline
$+$ & $\psdraw{a1}{0.4in}$ & $\psdraw{a2}{0.4in}$ & $\psdraw{a3}{0.4in}$ & 
$\psdraw{a4}{0.4in}$ & $\psdraw{a5}{0.4in}$ & $\psdraw{a6}{0.4in}$ \\ \hline
$R^+_v(s)$ & $q$ & $q-\bq$ & $1$ & $1$ & $0$ & $q$ \\ \hline
\end{tabular} }
\\
\text{\begin{tabular}{|c|c|c|c|c|c|c|} \hline
$-$ & $\psdraw{a1}{0.4in}$ & $\psdraw{a2}{0.4in}$ & $\psdraw{a3}{0.4in}$ & 
$\psdraw{a4}{0.4in}$ & $\psdraw{a5}{0.4in}$ & $\psdraw{a6}{0.4in}$ \\ \hline
$R^-_v(s)$ & $\bq $ & $0$ & $1$ & $1$ & $\bq-q$ & $ \bq $ \\ \hline
\end{tabular}
}
\end{gather}

where $\bq=q^{-1}$. (These weights come from the
$R$-matrix of the quantum group $U_q(\mathfrak{sl}_2)$ and we follow the
conventions of \cite[p.235]{CP}). 
\item
$\w(L)$ is the {\em writhe} of $D_L$, that is the sum of the signs of
the crossings of $D_L$.

\item  
To define the {\em rotation numbers} $\rot_0(s)$ and $\rot_1(s)$
of a state $s$, we first define the the rotation number $\rot(a)$ of a curve $a$ immersed in the plane to be the $\psi/(2\pi)$, where $\psi$ is the total rotation angle of the tangent vector of $a$. (The direction of the counter-clockwise rotation  is taken to be positive.) If $s$ is a state of $D_L$, then $\rot_i(s)$, for $i=1,2$,  is the sum of the rotation numbers of each of the $i$-colored circles in $s$. The rotation numbers can also be defined combinatorially as follows:
 smooth any 4-valent vertex
of $s$ as follows:
$$
\psdraw{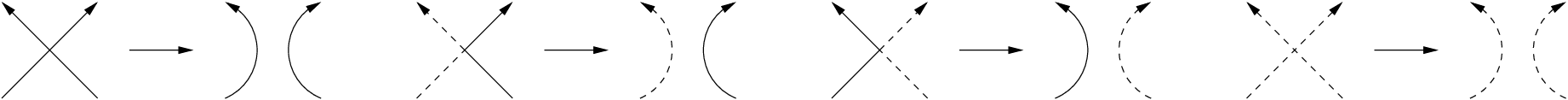}{5in}
$$
The result is a collection of oriented planar circles, colored by $0$ or $1$. 
$\rot_i(s)$ is the number of counter-clockwise $0$-colored circles minus
the number of clockwise $i$-colored circles. \end{itemize}

It will be convenient to include the phase factors $q^{\rot}$ in the $R$-matrix. This can be achieved as follows:  Let $\widetilde{D}_L$ be the diagram obtained by isotoping $D_L$ so that at each crossing the over and undercrossing arcs meet each other at a tangent where they cross.
This means that  the rotation
number of a cycle $c$ of $D_L$ equals  $\sum_{e \in c} \rot(\tilde{e})$, where $\rot(\tilde{e})$ is the rotation number of the edge $\tilde{e}$ of $\widetilde{D}_L$ corresponding to edge $e$ of $D_L$.

Now, for a state $s$ and a vertex  $v$ of $D_L$, consider its two
outgoing edges $e_1, e_2$ (see Figure \ref{fig:gadget}) 
and the rotation numbers $r_1=\rot(\tilde{e_1})$
and $r_2=\rot(\tilde{e_2})$. Of course, $r_1$ and $r_2$ depend on $v$. 
We now define the modified weights $B^{\pm}_v(s)$ of a state $s$ at a vertex
$v$ by:

\begin{gather}
\text{
\begin{tabular}{|c|c|c|c|c|c|c|} \hline
$+$ & $\psdraw{a1}{0.4in}$ & $\psdraw{a2}{0.4in}$ & $\psdraw{a3}{0.4in}$ & 
$\psdraw{a4}{0.4in}$ & $\psdraw{a5}{0.4in}$ & $\psdraw{a6}{0.4in}$ \\ \hline
$B^+_v(s)$ & $q q^{r_1+r_2}$ & $(q-\bq) q^{r_1-r_2}$ & $q^{-r_1+r_2}$ 
& $q^{r_1-r_2}$ & $0$ & $q q^{-r_1-r_2}$ \\ \hline
\end{tabular}
}\\
\text{
\begin{tabular}{|c|c|c|c|c|c|c|} \hline
$-$ & $\psdraw{a1}{0.4in}$ & $\psdraw{a2}{0.4in}$ & $\psdraw{a3}{0.4in}$ & 
$\psdraw{a4}{0.4in}$ & $\psdraw{a5}{0.4in}$ & $\psdraw{a6}{0.4in}$ \\ \hline
$B^-_v(s)$ & $\bq q^{r_1+r_2} $ & $0$ & $q^{-r_1+r_2} $ 
& $q^{r_1-r_2} $ & $(\bq-q) q^{-r_1+r_2}$ & $ \bq q^{-r_1-r_2}$ \\ \hline
\end{tabular}.
} \notag
\end{gather}

With these weights, we have
\begin{lemma}
\lbl{lem.b}
 Suppose that $D_L$ is a diagram for a link $L$. Let $\widetilde{D}_L$ be the diagram obtained by isotoping $D_L$ so that at each crossing the over and undercrossing arcs meet each other at a tangent where they cross. Then
$$
J(L)=J(\widetilde{D}_L)=q^{-2 \w(\widetilde{D}_L)} \sum_s \prod_v B^{\sgn(v)}_v(s).
$$
\end{lemma}

\begin{proof}
If $s$ is a state of $D_L$ and $s_i$ is the set of
$i$-colored circles in $s$, $i=0,1$ then 
$$
q^{\rot_0{s}-\rot_1(s)}=q^{\rot(s_0)-\rot(s_1)}=
\prod_{e} q^{\rot(\tilde{e}) \d_{e,s_0} - \rot(\tilde{e}) \d_{e,s_1}}, $$
where the second equality follows since the four arcs entering an exiting a gadget are tangenial. This sum can be written as
$$ \prod_v q^{\sum_{e \, \text{starts at} \, v} \rot(\tilde{e}) \d_{e,s_0} - \rot(\tilde{e}) \d_{e,s_1}}
$$
where $\d_{e,c}=1$ (respectively $0$) if $e$ lies in $c$ (respectively does not lie in $c$). 
Thus,
\begin{eqnarray*}
\sum_{s} 
q^{\rot_0(s)-\rot_1(s)} \prod_v R^{\sgn(v)}_v(s) &=&
\sum_s \prod_v q^{\sum_{e \, \text{starts at} \, v} \rot(\tilde{e}) \d_{e,s_0} - \rot(\tilde{e}) \d_{e,s_1}}
R^{\sgn(v)}_v(s) \\&=& \sum_s \prod_v B^{\sgn(v)}_v(s).
\end{eqnarray*}

\end{proof}

Note that in this proof we used the fact that in $\widetilde{D}_L$, the over and undercrossing arcs meet each other at a tangent where they cross.  We could avoid this tangent condition on the edges  entering and exiting a crossing by adding factors to the the expressions $B_v^{\pm}(s)$ that are determined by the angles formed by the crossings in $D_L$. See \cite{Ba} for details on this type of construction.

\begin{definition}
\lbl{def.ww}
Let $D_L$ be a link diagram. Then we let $\hat{D}_L$ denote the directed, edge-weighted graph drawn in the plane which is obtained from $D_L$ by replacing a neighbourhood of each crossing of $D_L$ with the gadgets and their incident edges as shown in Figure~\ref{fig:gadget}. The edge weights $W$ of $\hat{D}_L$ are given by the convention that all of  the undecorated edges have weight $1$,
\[ (a,b,c,d,e) = ((q-q^{-1}), (q+q^{-1})/2, -q/2, q, 1/2 ), \]

\[ (v,w,x,y,z)=(1/2,q^{-1}-q, (q+q^{-1})/2, -q^{-1}/2, q^{-1}), \]
and each edge $\hat{e}$ coming from an edge $e$ of $D_L$ is assigned the weight $q^{-2\rot (\tilde{e})}$. 
\end{definition}

\begin{lemma}
\lbl{lem.gadget}
Let $D_L$ be a link diagram of $L$ and $\hat{D}_L$ be the graph constructed as described in Definition~\ref{def.ww}. Then
\[ J(L)= q^{-2\omega (D_L)}  q^{\rot (L)} \sum_{p\in \mathcal{P}} \prod_{e\in p} W_e, \] 
where $\mathcal{P}$ denotes the set of vertex-disjoint directed cycles that cover all vertices of $\hat{D}_L$.
\end{lemma}
\begin{proof}
Let $s$ be a state of $D_L$ and $\tilde{s}$ be the corresponding state in $\tilde{D}_L$ (the diagram constructed in  Lemma~\ref{lem.b}). There is a correspondence between the state $\tilde{s}$ at a crossing $\tilde{v}$ and sets of vertex-disjoint directed cycles that cover all vertices of $\hat{D}_L$. This correspondence is given in the following way: if $\tilde{a}$ is an edge incident with $\tilde{v}$ in $\tilde{D}_L$, and  $\hat{a}$ is the corresponding edge in $\hat{D}_L$, then the edge $\hat{a}$ is in a cycle if and only if $\tilde{a}$ is colored by $1$.  The correspondence is given for a positive crossing in Appendix~\ref{s.f}.
From the appendix it is also readily seen that at a positive crossing 
\[  B_v^{+}(s) = \sum_{p\in \mathcal{P}(s,v)} q^{\rot(\hat{e}_1) +\rot(\hat{e}_2)) } \prod_{e\in p}W_e, \]
where if $\mathcal{P}(s)$ denotes the set of vertex-disjoint directed cycles that cover all vertices of $\hat{D}_L$ that corresponds with the state $s$, then  $\mathcal{P}(s,v)$  denotes the set of edges of  $\mathcal{P}(s)$ that belong to the gadget at the crossing $v$ of $D_L$.
It is readily checked that the corresponding identity also holds at negative crossings.
It then follows that 
\[ q^{-2\omega (D_L)}  q^{\rot (L)}  \sum_{p\in \mathcal{P}} \prod_{e\in p} W_e = q^{-2 \w(\widetilde{D}_L)} \sum_s \prod_v B^{\sgn(v)}_v(s),\]
which, by Lemma~\ref{lem.b}, is equal to the Jones polynomial, as required.
\end{proof}

We make the following definition for the matrix $M_L$ that is used in  Theorem~\ref{thm.per}.
\begin{definition}
Let $D_L$ be a link diagram and $\hat{D}_L$ be the associated graph constructed as in Definition~\ref{def.ww}. Then we let $M_L$  denote the adjacency matrix of $\hat{D}_L$, whose $(i,j)$-entry is the weight of the directed edge $(i,j)$.
\end{definition}

We can now prove Theorem~\ref{thm.per}.

\begin{proof}[Proof of Theorem~\ref{thm.per}.]
It is well known that if $A$ is the $0,1$
adjacency matrix of a graph, then $\Per(A)$ is the number of 
collections of vertex-disjoint directed cycles that cover all vertices of the graph. Similarly
for general matrix $A$, $\Per(A)=\sum_C a(C)$, where the sum is over all
collections of vertex-disjoint directed cycles that cover all vertices of the graph and
$a(C)=\prod_{\{ij\}\in C}A_{ij}$.  Using the notation of Lemma~\ref{lem.gadget}, it follows that  $\Per(M_L)=\sum_{p\in \mathcal{P}} \prod_{e\in p} W_e$. The result then follows  
from Lemma~\ref{lem.gadget}.
\end{proof}

\appendix

\section{Figures}
\lbl{s.f}

\bigskip

\begin{center}
\raisebox{7.5mm}{\epsfig{file=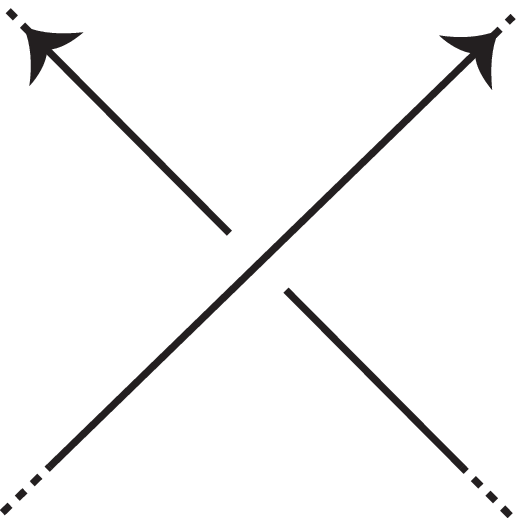, height=15mm} }  \raisebox{15mm}{:} \hspace{1cm}  \epsfig{file=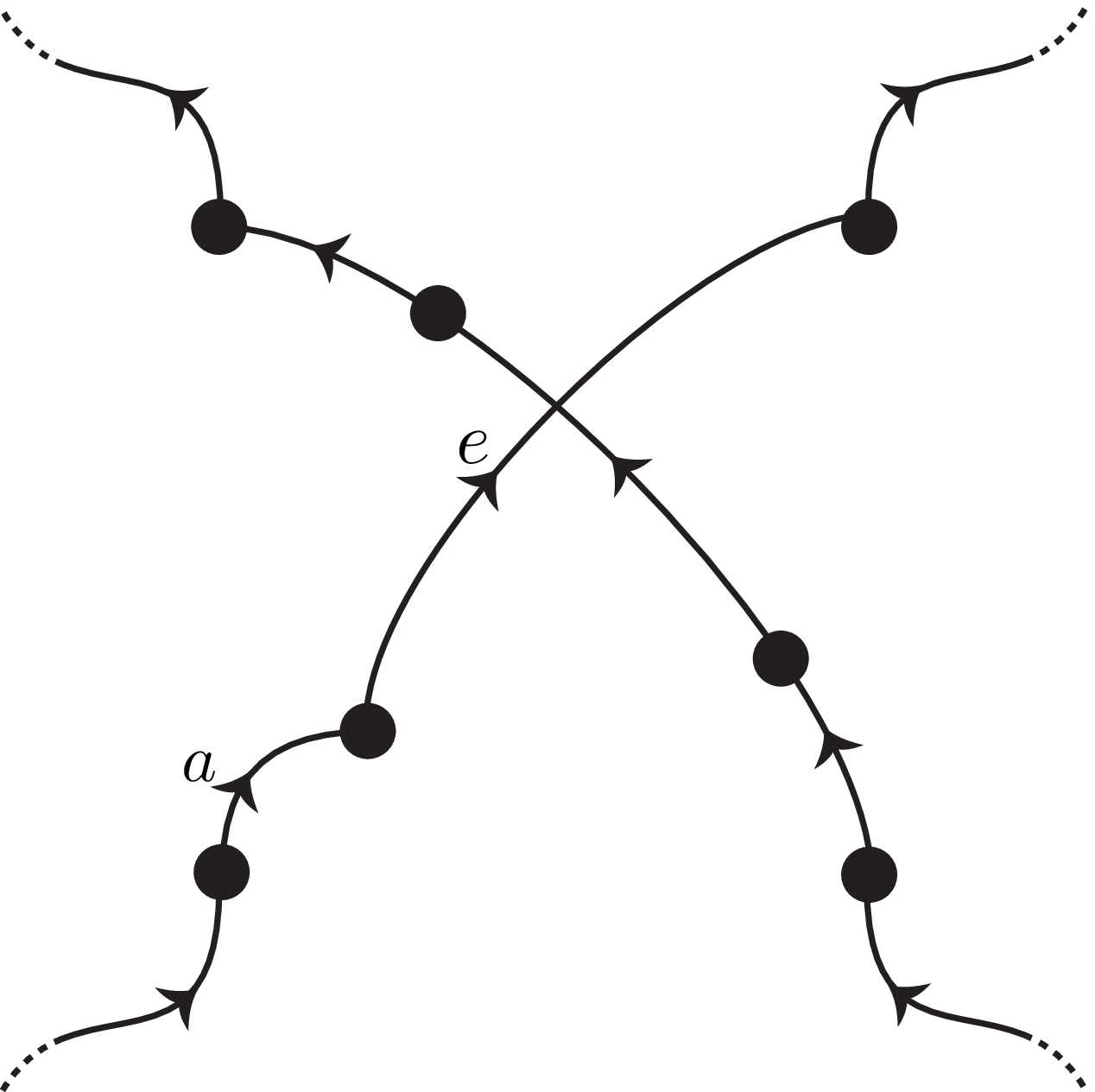, height=3cm} \hspace{5mm}  \raisebox{15mm}{+} \hspace{5mm} \epsfig{file=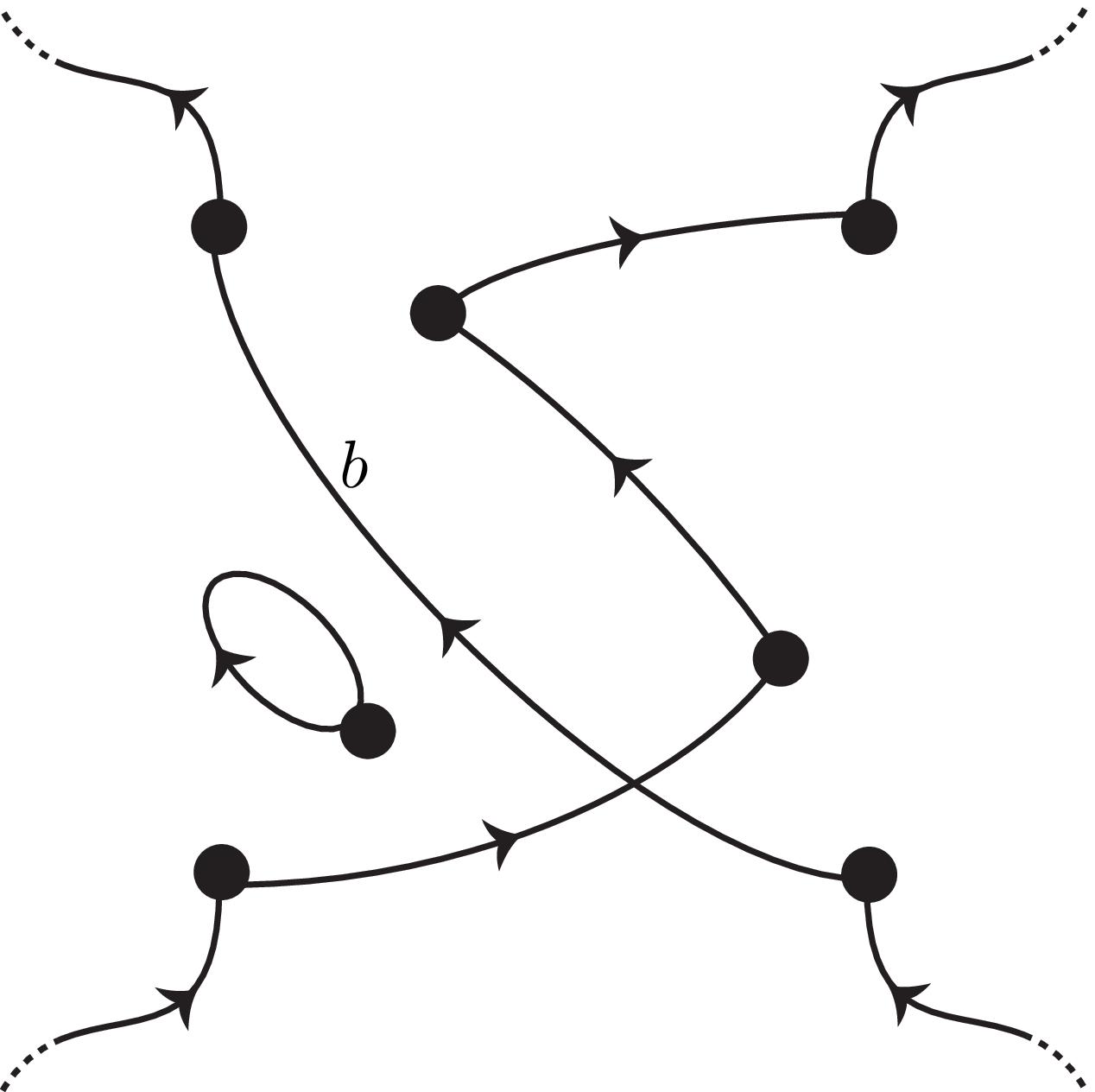, height=3cm} 
\end{center}
{\bf Case 1. } $q= ae+b$

\bigskip

\begin{center}
\raisebox{7.5mm}{\epsfig{file=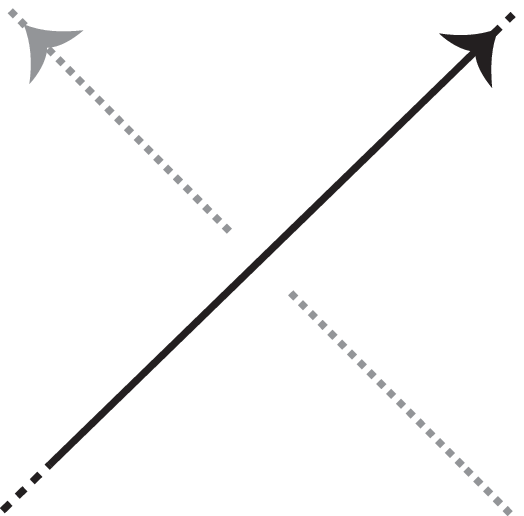, height=15mm} }  \raisebox{15mm}{:} \hspace{1cm}  \epsfig{file=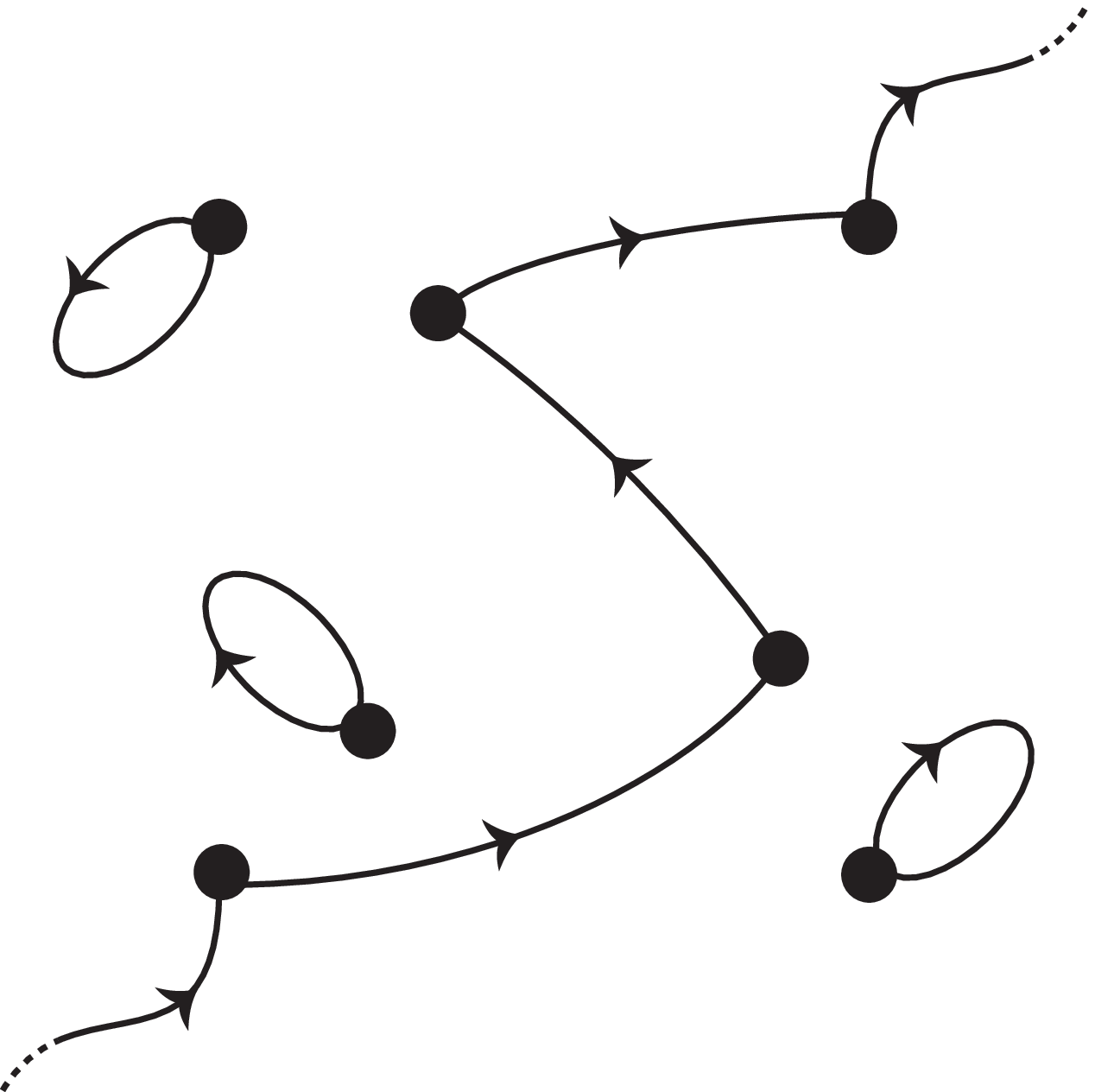, height=3cm} \end{center}
{\bf Case 2. } $1=1$
\bigskip

\begin{center}
\raisebox{7.5mm}{\epsfig{file=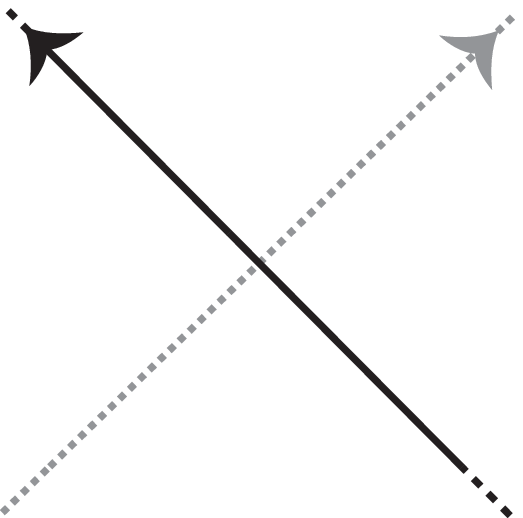, height=15mm} }  \raisebox{15mm}{:} \hspace{1cm}  \epsfig{file=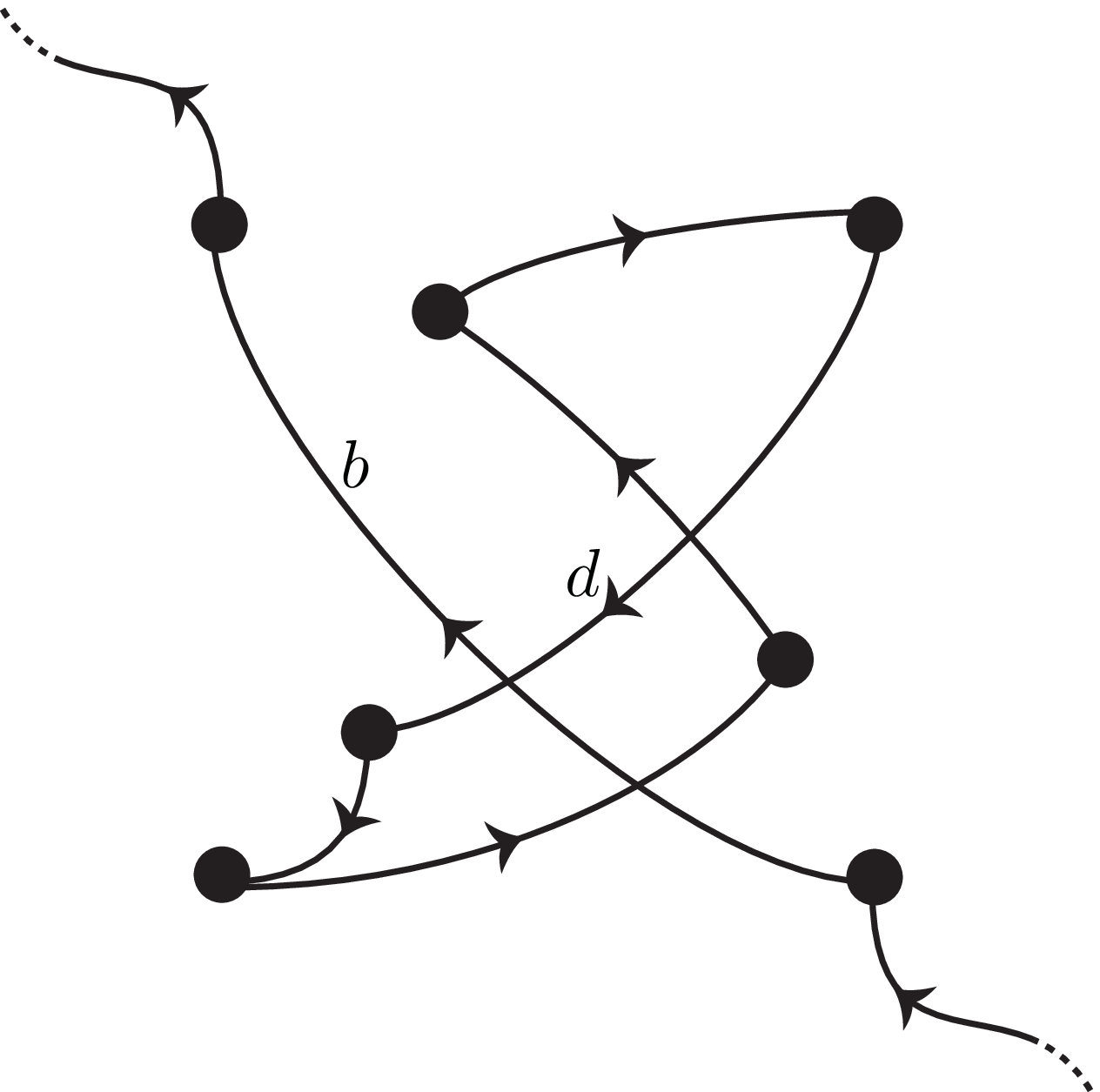, height=3cm} \hspace{5mm}  \raisebox{15mm}{+} \hspace{5mm} \epsfig{file=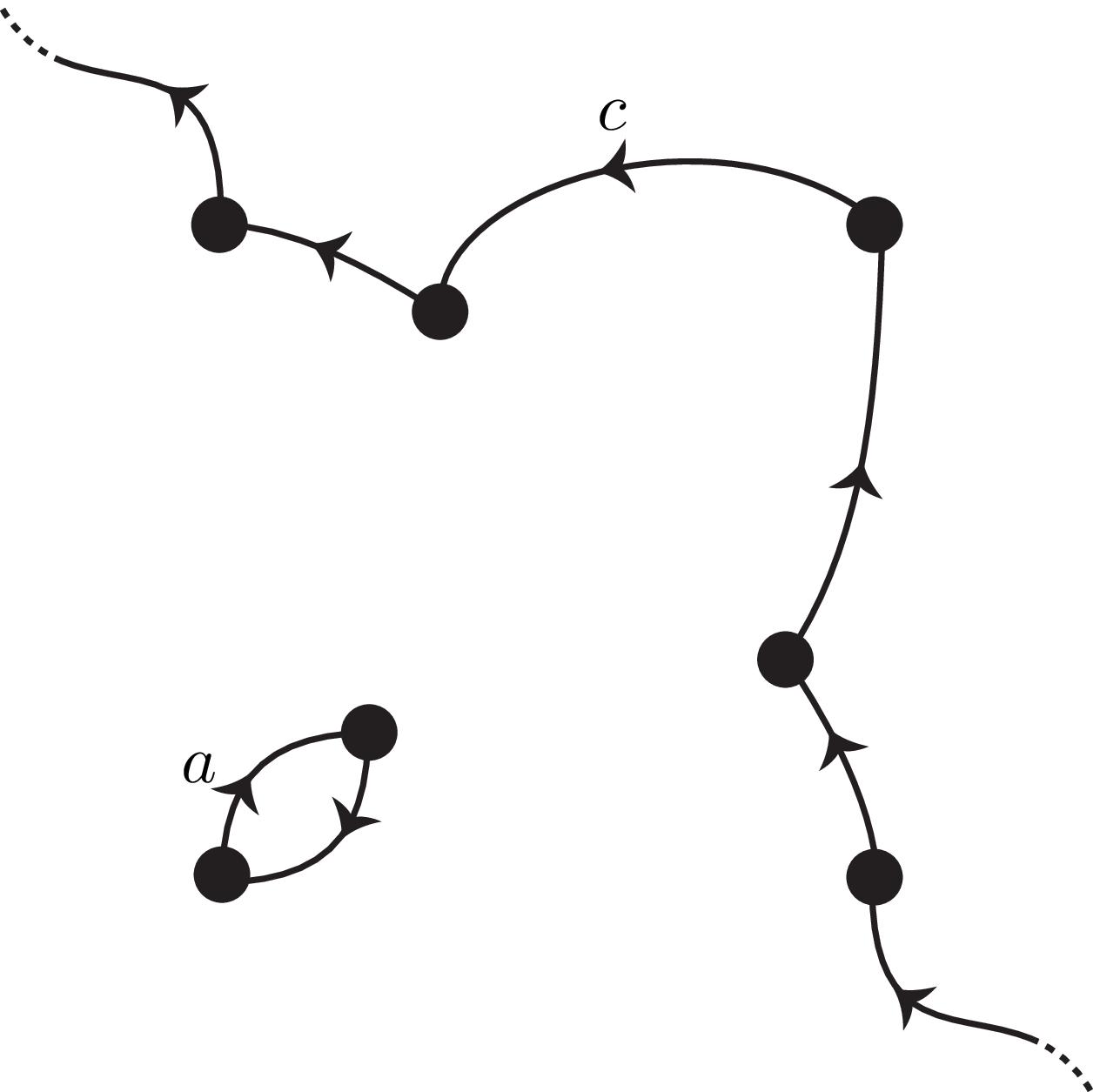, height=3cm} 
\end{center}
{\bf Case 3. } $1= bd+ ac$
\bigskip

\begin{center}
\raisebox{6mm}{\epsfig{file=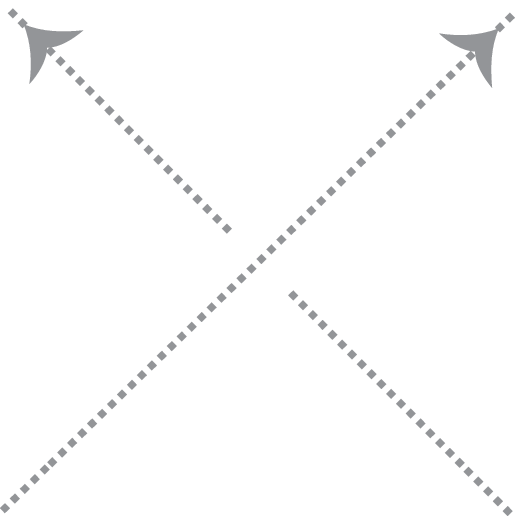, height=15mm} }  \raisebox{12mm}{:} \hspace{1cm}   \epsfig{file=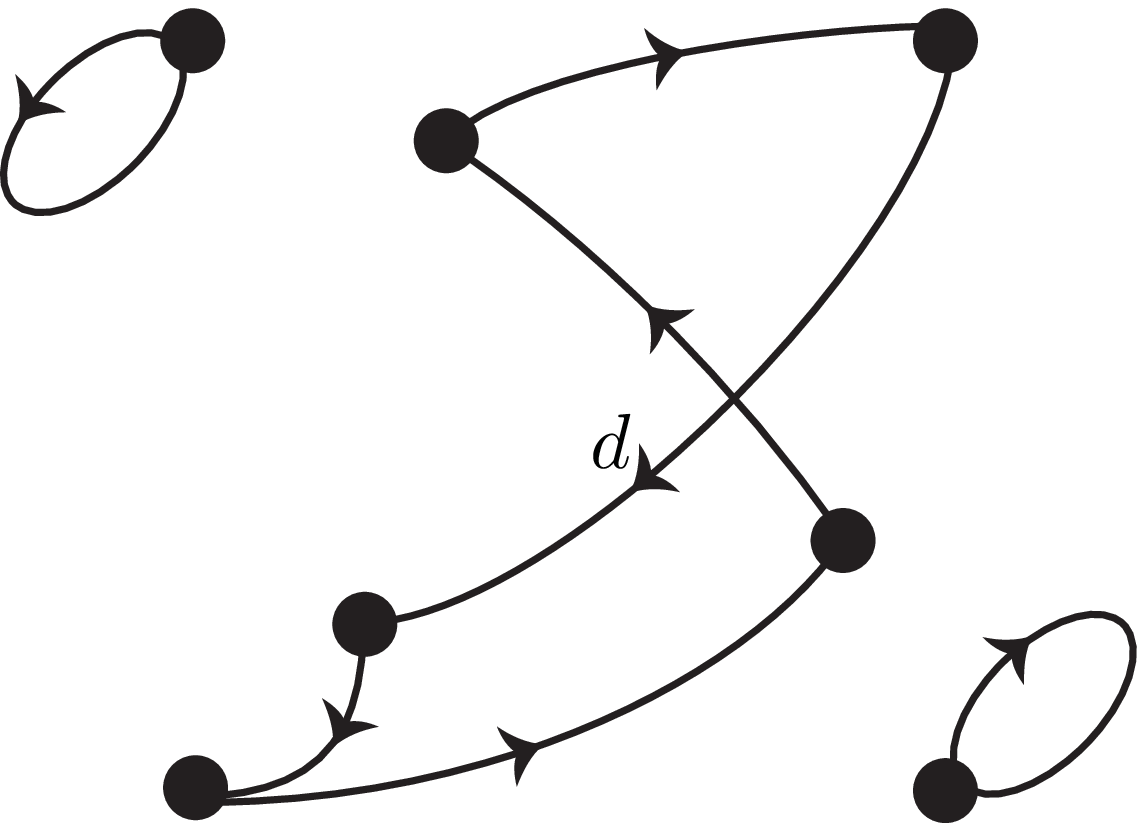, height=2.4cm} 
\end{center}
{\bf Case 4. } $q=d$
\bigskip

\begin{center}
\raisebox{7.5mm}{\epsfig{file=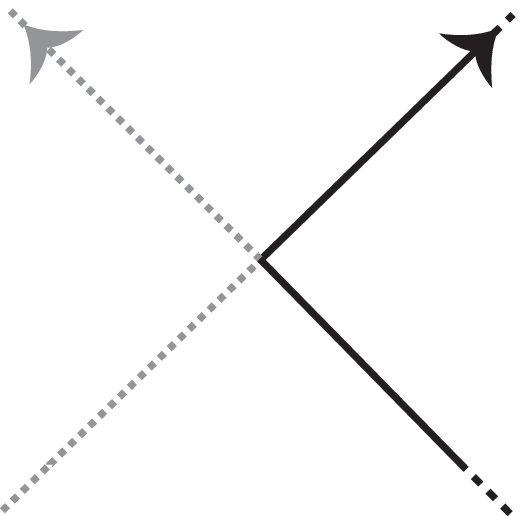, height=15mm} }  \raisebox{15mm}{:} \hspace{1cm}  \epsfig{file=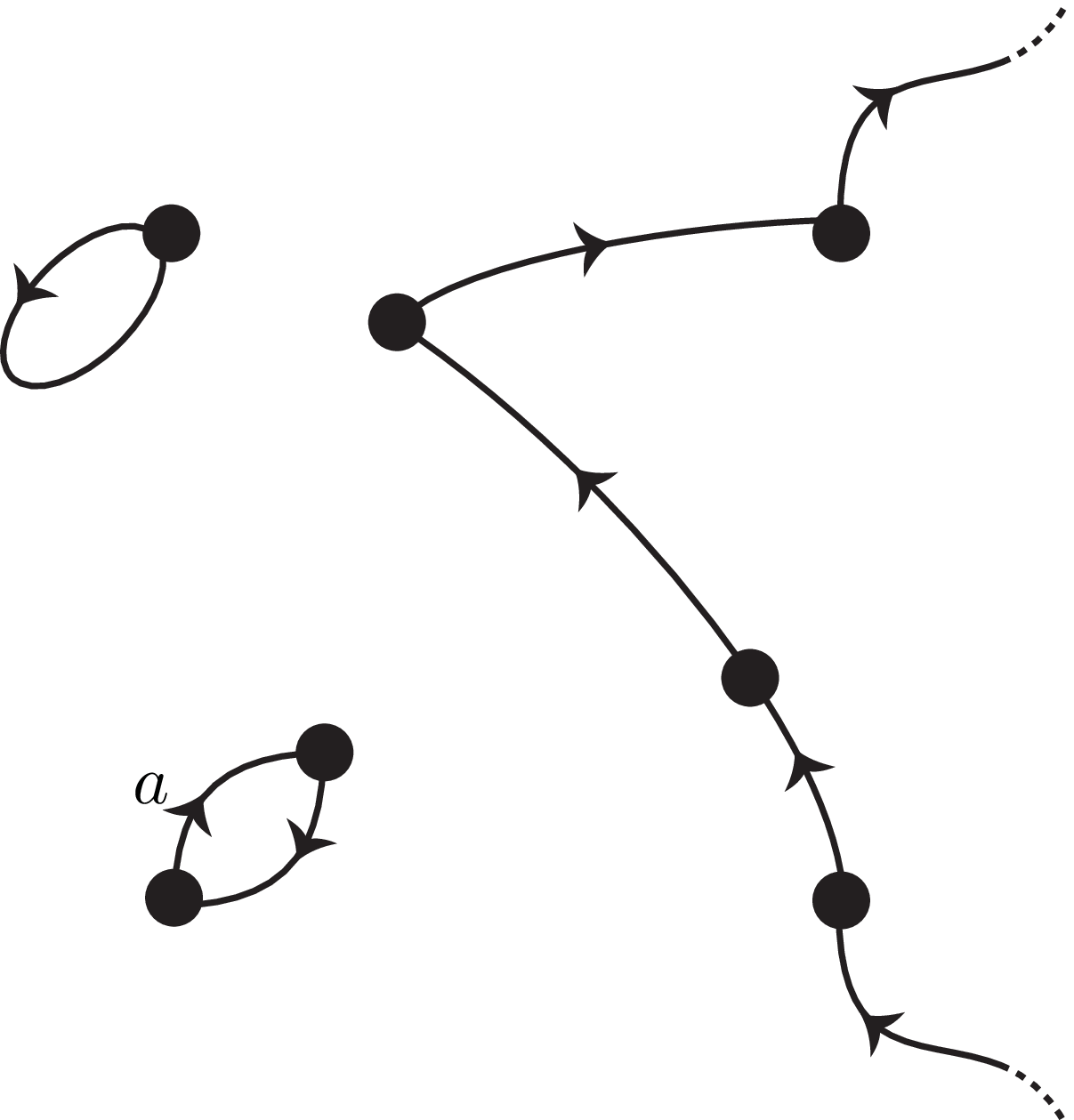, height=3cm} \hspace{5mm} 
\end{center}
{\bf Case 5.} $q-q^{-1}= a$
\bigskip

\begin{center}
\raisebox{7.5mm}{\epsfig{file=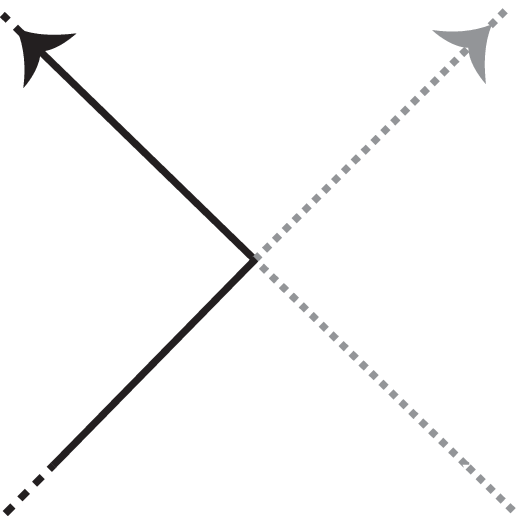, height=15mm} }  \raisebox{15mm}{:} \hspace{1cm}  \epsfig{file=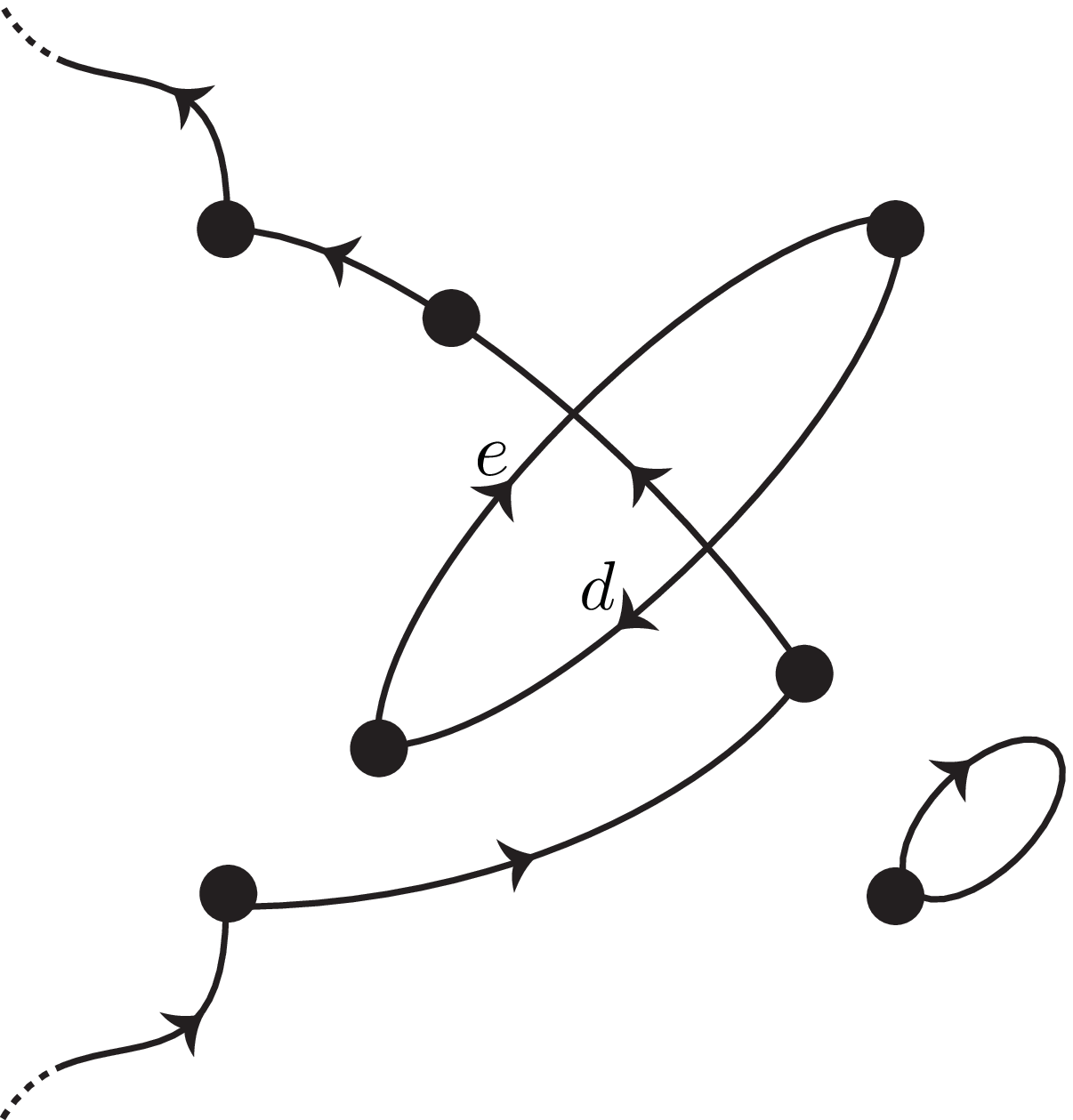, height=3cm} \hspace{5mm}  \raisebox{15mm}{+} \hspace{5mm} \epsfig{file=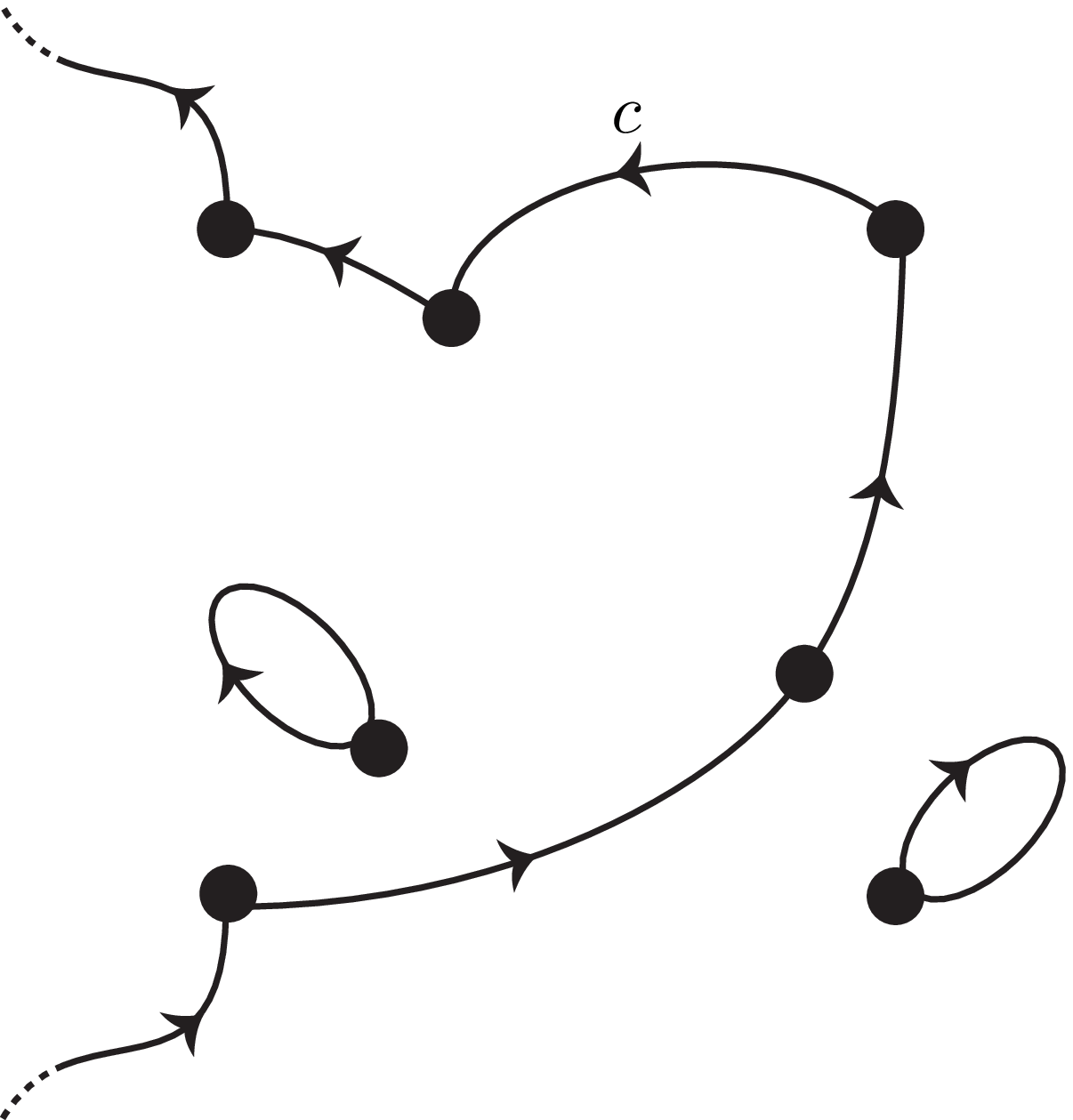, height=3cm} 
\end{center}
{\bf Case 6.} $0= de+ c$

%






\ifx\undefined\bysame
	\newcommand{\bysame}{\leavevmode\hbox
to3em{\hrulefill}\,}
\fi

\end{document}